\documentclass[12pt]{amsart}
\usepackage{amssymb}
\usepackage{mathrsfs}
\usepackage[all]{xy}

\newcommand{\Weiss}{
\author{Tomasz Weiss}
\address{Institute of Mathematics, Akademia Podlaska
08-119 Siedlce, Poland}
\email{tomaszweiss@go2.pl}
}

\newcommand{\Tsaban}{
\author{Boaz Tsaban}
\thanks{Supported by the Koshland Center for Basic Research.}
\address{Department of Mathematics,
Weizmann Institute of Science, Rehovot 76100, Israel}
\email{boaz.tsaban@weizmann.ac.il}
\urladdr{http://www.cs.biu.ac.il/\~{}tsaban} }

\newcommand{\Binary}{T}
\newcommand{\gpbl}{{\mbox{\textit{\tiny gp}}}}
\newcommand{\two}{{\{0,1\}}}

\newcommand{\productive}{semiproductive}
\renewcommand{\|}{\rest} 
\newcommand{\Pgal}{{\tilde P}}
\newcommand{\isoproductive}{iso-productive}
\newcommand{\zeroproductive}{$0$-productive}

\newcommand{\HURE}{\mathcal{H}}

\newcommand{\SMZ}{\mathrm{SMZ}}

\newcommand{\bq}{\begin{quote}}
\newcommand{\eq}{\end{quote}}
\newcommand{\inv}{^{-1}}

\long\def\forget#1\forgotten{}

\newcommand{\diam}{\op{diam}}
\newcommand{\cC}{\mathcal{C}}
\newcommand{\cU}{\mathcal{U}}

\newcommand{\cP}{\mathcal{P}}

\newcommand{\seq}[1]{\{#1\}_{n\in\N}}

\newcommand{\op}{\operatorname}

\newcommand{\Cantor}{{\two^\N}}

\newcommand{\Cal}{\mathcal}
\newcommand{\Madd}{{\cM}^\star}

\newcommand{\scrA}{\mathscr{A}}
\newcommand{\scrB}{\mathscr{B}}

\newcommand{\cF}{{\Cal F}}

\newcommand{\cM}{\mathcal{M}}
\newcommand{\cN}{\mathcal{N}}
\newcommand{\N}{\mathbb{N}}

\newcommand{\Null}{\mathcal{N}}
\newcommand{\NN}{{\N^\N}}
\newcommand{\NNup}{\N^{\nearrow\N}}

\renewcommand{\O}{\Cal O}

\newcommand{\Q}{\mathbb{Q}}
\newcommand{\R}{\mathbb{R}}

\newcommand{\Union}{\bigcup}
\newcommand{\Z}{\mathbb{Z}}
\newcommand{\Impl}{\Rightarrow}

\renewcommand{\b}{{\mathfrak b}}

\renewcommand{\c}{{\mathfrak c}}

\renewcommand{\i}{\item}

\newcommand{\oo}{\infty}

\newcommand{\x}{\times}

\newcommand{\nin}{\not\in}

\newcommand{\sbst}{\subseteq}

\newcommand{\sm}{\setminus}

\renewcommand{\(}{\left(}
\renewcommand{\)}{\right)}
\newcommand{\<}{\langle}
\renewcommand{\>}{\rangle}

\newcommand{\cI}{\mathcal{I}}
\newcommand{\cJ}{\mathcal{J}}
\newcommand{\cK}{\mathcal{K}}

\newcommand{\cov}{{\sf cov}}

\newtheorem{thm}{Theorem}[section]

\newtheorem{prop}[thm]{Proposition}
\newtheorem{prob}[thm]{Problem}

\newtheorem{lem}[thm]{Lemma}

\theoremstyle{definition}
\newtheorem{definition}[thm]{Definition}
\newtheorem{abbr}[thm]{Abbreviation}

\theoremstyle{remark}
\newtheorem{rem}[thm]{Remark}

\newcommand{\be}{\begin{enumerate}}
\newcommand{\ee}{\end{enumerate}}
\newcommand{\bi}{\begin{itemize}}
\newcommand{\ei}{\end{itemize}}

\newcommand{\sone}{{\sf S}_1}    \newcommand{\sfin}{{\sf S}_{fin}}
\newcommand{\ufin}{{\sf U}_{fin}}

\newcommand{\closure}{\overline}

\newcommand{\rest}{{\mathord{\restriction}}}
\newcommand{\lft}[2]{\mathopen\ifcase#1{}\oo\or
                        \big#2\or\Big#2\else\oo\fi}
\newcommand{\rgt}[2]{\mathclose\ifcase#1{}\oo\or
                        \big#2\or\Big#2\else\oo\fi}

\title{Products of special sets of real numbers}

\Tsaban
\Weiss



\begin{document}
\begin{abstract}
We describe a simple machinery which translates results on algebraic
sums of sets of reals into the corresponding results on their
cartesian product. Some consequences are:
\be
\i The product of a meager/null-additive set
and a strong measure zero/strongly meager set in the Cantor space
has strong measure zero/is strongly meager, respectively.
\i Using Scheepers' notation for selection principles:
$\sfin(\Omega,\Omega^\gpbl)\cap\sone(\O,\O)=\sone(\Omega,\Omega^\gpbl)$, and
Borel's Conjecture for $\sone(\Omega,\Omega)$ (or just $\sone(\Omega,\Omega^\gpbl)$)
implies Borel's Conjecture.
\ee
These results extend results of Scheepers and Miller, respectively.
\end{abstract}

\maketitle

\section{Products in the Cantor space}\label{cantorprods}

The Cantor space $\cC=\Cantor$ is equipped with the product topology.
For distinct $x,y\in \cC$, write $N(x,y)=\min\{n : x(n)\neq y(n)\}$.
Then the topology of $\cC$ is generated by the following metric:
$$d(x,y) =
\begin{cases}
\frac{1}{N(x,y)+1} & x \neq y\\
0 & x=y
\end{cases}$$
(so that $d(x,y)\le 1$ for all $x,y\in\cC$). A canonical measure
$\mu$ is defined on $\cC$ by taking the product of the uniform
probability measure on $\two$. Fix a natural number $k$, and
consider the product space $\cC^k$. Define the \emph{product
metric} $d_k$ on $\cC^k$ by
$$d_k((x_0,\dots,x_{k-1}),(y_0,\dots,y_{k-1})) = \max\{d(x_0,y_0),\dots,d(x_{k-1},y_{k-1})\}.$$
Then $d_k$ generates the topology of $\cC^k$.
The measure on $\cC^k$ is the product measure $\mu\x \dots\x \mu$ ($k$ times).

$\cC$, with the operation $\oplus$
defined by $(x\oplus y)(n) = x(n)+y(n)\bmod 2$
is a topological group, and therefore so is
$\cC^k$ for all $k$.
\begin{lem}
The function $\Psi_k:\cC^k\to\cC$ defined by
$$\Psi_k(x_0,\dots,x_{k-1})(mk+i) = x_i(m)$$
for each $m$ and each $i<k$,
is a bi-Lipschitz measure preserving group isomorphism.
\end{lem}
\begin{proof}
Clearly $\Psi_k$ is bijective.
Assume that $\vec x = (x_0,\dots,x_{k-1})$ and
$\vec y = (y_0,\dots,y_{k-1})$ are members of $\cC^k$.
Then for each $m$ and each $i<k$,
\begin{eqnarray*}
\lefteqn{\Psi_k(\vec x \oplus \vec y)(mk+i) = }\\
 & = & \Psi_k(x_0 \oplus y_0,\dots,x_{k-1}\oplus y_{k-1})(mk+i) = (x_i \oplus y_i)(m) =\\
 & = & x_i(m) + y_i(m) \bmod 2 = \Psi_k(\vec x)(mk+i) + \Psi_k(\vec y)(mk+i) \bmod 2.
\end{eqnarray*}
Thus, $\Psi_k(\vec x \oplus \vec y) = \Psi_k(\vec x) \oplus \Psi_k(\vec y)$,
and $\Psi_k$ is a group isomorphism.

Now, assume that $\vec x = (x_0,\dots,x_{k-1})$ and
$\vec y = (y_0,\dots,y_{k-1})$ are distinct members of $\cC^k$,
and let $i$ be such that $d(x_i,y_i)$ is maximal, that is,
$N=N(x_i,y_i)$ is minimal.
Then $N(\Psi_k(\vec x),\Psi_k(\vec y)) \ge kN$,
and therefore
$$d(\Psi_k(\vec x),\Psi_k(\vec y))\le \frac{1}{kN+1} < d(x,y).$$
Similarly, for distinct $x,y\in\cC$, if $N(x,y) = mk+i$ where
$i<k$, then $N(\Psi_k\inv(x),\allowbreak\Psi_k\inv(y))\ge m$ and
$$d(\Psi_k\inv(x),\Psi_k\inv(y))\le \frac{1}{m+1} \le \frac{k}{mk+i+1} = k\cdot d(x,y).$$

To see that $\Psi_k$ is measure preserving, observe that the measure of a basic open set $U$
in $\cC$ is $2^{-m}$, where $m$ is the number of coordinates of $U$ which are not equal to $\two$.
Consequently, the same assertion is true for $\cC^k$, where $m$ is the sum of numbers of such coordinates
within each of the $k$ coordinates of $\cC^k$. This number $m$
is invariant under $\Psi_k$; thus $\Psi_k$ preserves measures of basic open sets.
\end{proof}
This often allows us to restrict attention to subsets of $\cC$
rather than subsets of $\cC^k$ for arbitrary $k$.


\begin{abbr}
$\Pgal(A)=\Union_{k\in\N}P(A^k)$.
\end{abbr}

\begin{thm}\label{prod2sum}
Assume that $\cI,\cJ,\cK\sbst\Pgal(\cC)$, and that $\cK$ is closed under taking
Lipschitz images.
If
\bq
for each $X\in\cI\cap P(\cC^k)$ and $Y\in\cJ\cap P(\cC^k)$, $X\x Y\in\cK$,
\eq
then
\bq
for each $X\in\cI\cap P(\cC^k)$ and $Y\in\cJ\cap P(\cC^k)$, $X\oplus Y\in\cK$.
\eq
\end{thm}
\begin{proof}
The mapping $(x,y)\mapsto x\oplus y$ is a Lipschitz mapping
from $X\x Y$ onto $X\oplus Y$.
\end{proof}

The converse of Theorem \ref{prod2sum} also holds, and in a much stronger form.
For simplicity, we introduce the following notions.
\begin{abbr}
$\cP\sbst\Pgal(\cC)$
is \emph{\productive{}} if:
\be
\i For each $k$, $l$, and $X\in\cP\cap\cC^k$,
   if $0$ is the zero element of $\cC^l$, then $X\x\{0\}\in\cP$; and
\i For each $k$, $l$, $X\in\cP\cap\cC^k$, and a
   bi-Lipschitz measure preserving group isomorphism $\Phi:\cC^k\to\cC^l$,
   $\Phi[X]\in\cP$.
\ee
We will say that $\cP$ is \emph{\zeroproductive} if we only require that (1) is satisfied,
and \emph{\isoproductive{}} if we only require that (2) is satisfied.
\end{abbr}

Many properties of
special sets of reals are \productive{}, e.g.,
Hausdorff dimension, strong measure zero, the properties in
Cicho\'n's Diagram for small sets \cite{pawlikowskireclaw}
or in Scheepers' Diagram (see Section \ref{BC}) and its extensions
\cite{CBC, tautau}; see \cite{MilSpec, coc1} for more examples.

As changing the order of coordinates is a bi-Lipschitz
measure preserving group isomorphism,
we have the following.
\begin{lem}
Assume that $\cP\sbst\Pgal(\cC)$ is \productive{}.
Then for each $k$, $l$, and $X\in\cP\cap\cC^k$,
if $0$ is the zero element of $\cC^l$, then $\{0\}\x X\in\cP$.
\hfill\qed
\end{lem}

\begin{thm}\label{sum2prod}
Assume that $\cI, \cJ, \cK\sbst\Pgal(\cC)$ such that $\cI$ and $\cJ$
are \productive{}, and $\cK$ is \isoproductive{}.
If
\bq
for each $X\in\cI\cap P(\cC)$ and $Y\in\cJ\cap P(\cC)$, $X\oplus Y\in\cK$,
\eq
then
\bq
for each $X\in\cI$ and $Y\in\cJ$, $X\x Y\in\cK$.
\eq
\end{thm}
\begin{proof}
Assume that $X\in\cI\cap P(\cC^k)$ and $Y\in\cJ\cap P(\cC^l)$.
Then $\tilde X = \Psi_k[X]\in\cI\cap P(\cC)$, and $\tilde Y = \Psi_l[Y]\in\cJ\cap P(\cC)$.
\begin{eqnarray*}
\Psi_2[\tilde X\x \tilde Y] & = & \Psi_2[(\tilde X\x\{0\})\oplus(\{0\}\x \tilde Y)] =\\
 & = & \Psi_2[\tilde X\x\{0\}]\oplus\Psi_2[\{0\}\x \tilde Y]= X'\oplus Y',
\end{eqnarray*}
Thus, $X'\oplus Y'\in\cK$. As $\Psi_2$ is bijective,
$\tilde X\x \tilde Y = \Psi_2\inv[X'\oplus Y']\in\cK$.
As $\Psi_k\x \Psi_l$ is a bi-Lipschitz measure preserving group isomorphism and
$\Psi_k\inv\x\Psi_l\inv:\tilde X\x\tilde Y\to X\x Y$ is surjective,
$X\x Y\in\cK$.
\end{proof}

Assume that $\cI\sbst\Pgal(\cC)$.
A subset $X$ of $\cC^k$ is \emph{$\cI$-additive} if for each $I\in\cI\cap P(\cC^k)$,
$X\oplus I\in\cI$. Clearly if $X,Y\sbst\cC$ are $\cI$-additive, then
$X\oplus Y$ is $\cI$-additive.
More generally,
a subset $X$ of $\cC^k$ is \emph{$(\cI,\cJ)$-additive} if for each $I\in\cI\cap\cC^k$,
$X\oplus I\in\cJ$.
Let $\cI^\star$ and $(\cI,\cJ)^\star$ denote the classes of all $\cI$-additive sets and
$(\cI,\cJ)$-additive sets, respectively.

\begin{thm}\label{prods*}
Assume that $\cI,\cJ,\cK\sbst\Pgal(\cC)$ are \isoproductive{},
and $(\cJ,\cK)^\star, (\cI,\cJ)^\star$ are \zeroproductive{}.
Then for each $X\in(\cJ,\cK)^\star$ and $Y\in(\cI,\cJ)^\star$, $X\x Y\in (\cI,\cK)^\star$.
In particular, if $\cI$ is \isoproductive{} and $\cI^\star$ is \zeroproductive{}, then
$\cI^\star$ is closed under taking finite products.
\end{thm}
\begin{proof}
Assume that $X\in(\cJ,\cK)^\star\cap P(\cC)$, $Y\in(\cI,\cJ)^\star\cap P(\cC)$, and $I\in\cI$.
Then $Y\oplus I\in\cJ$ and therefore $X\oplus (Y\oplus I)\in\cK$.
Thus $X\oplus Y\in(\cI,\cK)^\star$.

\begin{lem}\label{iso*iso}
Assume that $\cI,\cJ\sbst\Pgal(\cC)$ are \isoproductive{}.
Then $(\cI,\cJ)^\star$ is \isoproductive{}.
\end{lem}
\begin{proof}
Assume that $X\in(\cI,\cJ)^\star\cap P(\cC^k)$ and
$\Phi:\cC^k\to \cC^l$ is a bi-Lipschitz measure preserving group isomorphism.
Then for each $I\in\cI\cap P(\cC^l)$,
$$\Phi[X]\oplus I = \Phi[X\oplus\Phi\inv[I]].$$
As $\Phi\inv[I]\in\cI$, $X\oplus\Phi\inv[I]\in\cJ$ and
therefore $\Phi[X]\oplus I\in\cJ$, that is,
$\Phi[X]\in (\cI,\cJ)^\star$.
\end{proof}
Thus $(\cJ,\cK)^\star$ and $(\cI,\cJ)^\star$ are \productive{},
and the theorem follows from Theorem \ref{sum2prod}.
\end{proof}

Following is a useful criterion for the $0$-productivity
required in Theorem \ref{prods*}.

\begin{lem}\label{0prods*}
Assume that
$\cI,\cJ\sbst\Pgal(\cC)$ are \isoproductive{},
and for the element $0\in\cC$ and
each $X\in(\cI,\cJ)^\star\cap P(\cC)$,
$\Psi_2[X\x\{0\}]\in(\cI,\cJ)^\star$.
Then $(\cI,\cJ)^\star$ is \zeroproductive{}.
In particular, if $\cI$ is \isoproductive{} and
and for each $X\in\cI^\star\cap P(\cC)$,
$\Psi_2[X\x\{0\}]\in\cI^\star$,
then $\cI^\star$ is \zeroproductive{}.
\end{lem}
\begin{proof}
Assume that $X\in(\cI,\cJ)^\star\cap P(\cC^k)$ and fix $l$.
By Lemma \ref{iso*iso}, $\tilde X=\Psi_k[X]\in(\cI,\cJ)^\star\cap P(\cC)$.
Thus $\Psi_2\circ(\Psi_k\x\Psi_l)[X\x\{0\}]=\Psi_2[\tilde X\x\{0\}]\in(\cI,\cJ)^\star$.
Applying Lemma \ref{iso*iso} again, we get that
$X\x\{0\}\in(\cI,\cJ)^\star$.
\end{proof}

We now give some applications.
Let $X$ be a metric space.
Following Borel,
we say that $X$ has \emph{strong measure zero} if
for each sequence of positive reals $\seq{\epsilon_n}$,
there exists a cover
$\seq{I_n}$ of $X$ such that $\diam(I_n)<\epsilon_n$ for all $n$.
$X$ has the \emph{Hurewicz property}
if for each sequence $\seq{\cU_n}$ of open covers of $X$
there exist finite subsets $\cF_n\sbst\cU_n$, $n\in\N$, such
that $X\sbst\bigcup_n\bigcap_{m>n} \cup\cF_n$.
Let $\SMZ$ (respectively, $\HURE$) denote the collections of metric spaces having
strong measure zero (respectively, the Hurewicz property).

The following theorem of Scheepers will serve as a ``test case'' for
our approach.

\begin{thm}[Scheepers \cite{smzpow}]\label{schthm}
Let $X$ be a strong measure zero metric
space which also has the Hurewicz property. Then for each strong
measure zero metric space $Y$, $X\times Y$ has strong measure zero.
\end{thm}
Scheepers' proof of Theorem \ref{schthm} is by a reduction of the
Hurewicz property to the so called ``grouping property'', which is
proved using a result from topological game-theory. We will
present alternative proofs for the case that the spaces are sets
of real numbers. We first observe that this follows from the
corresponding theorem with $X\oplus Y$ instead of $X\x Y$, which was
proved in \cite{NSW, pawlikowskireclaw}:
Since the collections $\HURE\cap\Pgal(\cC)$
and $\SMZ\cap\Pgal(\cC)$ are \productive{}, Theorem \ref{sum2prod}
applies. (For another proof of Scheepers' Theorem in $\cC$, see
Theorem \ref{barthm} in the appendix.)

We now treat the classes of meager-additive and null-additive sets.
Let $\cM$ and $\cN$ denote the meager (i.e., first category) and null
(i.e., measure zero) sets, respectively.
\begin{thm}\label{M&Nstarprod}
$\Madd \cap\Pgal(\cC)$ and $\cN^\star\cap\Pgal(\cC)$ are \productive{}, and are
closed under taking finite products.
\end{thm}
\begin{proof}
Clearly, $\cM\cap\Pgal(\cC)$ and $\cN\cap\Pgal(\cC)$ are \productive{}.
By Theorem \ref{prods*}, it is enough to show that
the classes $\Madd \cap\Pgal(\cC)$ and $\cN^\star\cap\Pgal(\cC)$ are \zeroproductive{}.
We first treat $\Madd $.

\begin{lem}[{Bartoszy\'nski-Judah \cite[Theorem 2.7.17]{jubar}}]\label{M*char}
A subset $X$ of $\cC$ is meager-additive if, and only if,
for each increasing sequence $\seq{m_n}$ there exist
a sequence $\seq{l_n}$ and $y\in\cC$ such that for each
$x\in X$ and all but finitely many $n$,
$$l_n\le m_k<m_{k+1}\le l_{n+1}\mbox{ and } x\|[m_k,m_{k+1})=y\|[m_k,m_{k+1})$$
for some $k$.
(In this case we say that $\seq{l_n}$ and $y$ are appropriate for $\seq{m_n}$ and $X$.)
\hfill\qed
\end{lem}

We will prove the sufficient criterion of Lemma \ref{0prods*}.
Assume that $X\in\Madd \cap P(\cC)$, and let
$\tilde X=\Psi_2[X\x\{0\}]$. We must show that $\tilde X\in\Madd $.
Let an increasing sequence $\seq{m_n}$ be given.
Choose an increasing sequence $\seq{m_n'}$ of even numbers
such that for all but finitely many $n$,
there exists $k$ such that $m_n'\le m_k<m_{k+1}\le m_{n+1}'$.

Apply Lemma \ref{M*char} to obtain $\seq{l_n}$ and $y$ which are
appropriate for $\seq{m_n'/2}$ and $X$.
By the definition of $\Psi_2$, we get that $\seq{2l_n}$ and $\Psi_2(y,0)$
are appropriate for $\seq{m_n'}$ and $\Psi_2[X\x\{0\}]$.
In particular, they are appropriate for $\seq{m_n}$ and $\Psi_2[X\x\{0\}]$.
This shows that $\Madd $ is \zeroproductive{}.

Using similar arguments, the fact that $\cN^\star\cap\Pgal(\cC)$ is \zeroproductive{} follows
from the following.
\begin{lem}[{Shelah \cite[Theorem 2.7.18]{jubar}}]
A subset $X$ of $\cC$ is null-additive if, and only if,
for every increasing sequence $\seq{m_n}$ there exists a sequence
$\seq{S_n}$ such that each $S_n$ is a set of at most $n$
functions from $[m_n,m_{n+1})$ to $\two$, and for each
$x\in X$ and all but finitely many $n$, $x\|[m_n,m_{n+1})\in S_n$.
\hfill\qed
\end{lem}
This finishes the proof of Theorem \ref{M&Nstarprod}.
\end{proof}

\begin{prop}[folklore]\label{GMSk}
For all $k$, a set $X\sbst\cC^k$ has strong measure zero if, and only if,
for each meager $M\sbst\cC^k$, $X\oplus M\neq\cC^k$.
\end{prop}
\begin{proof}
Assume that $X\sbst\cC^k$ has strong measure, and $M\sbst\cC^k$
is meager.
Then by the Galvin-Mycielski-Solovay Theorem, 
$\Psi_k[X\oplus M] = \Psi_k[X]\oplus\Psi_k[M]\neq\cC$,
and therefore $X\oplus M\neq\cC^k$.

Conversely, assume that $X\sbst\cC^k$ and
for each meager $M\sbst\cC^k$, $X\oplus M\neq\cC^k$.
Then for each meager $M\sbst\cC$,
$X\oplus\Psi_k\inv[M]\neq\cC^k$, therefore
$\Psi_k[X\oplus\Psi_k\inv[M]] = \Psi_k[X]\oplus M\neq\cC$,
thus $\Psi_k[X]$ has strong measure zero, and therefore
$X$ has strong measure zero.
\end{proof}

Every set of reals with the Hurewicz property
as well as strong measure zero is meager-additive (\cite{NSW}, or
Theorem \ref{(*)impliesM*} below).
Consequently, the following theorem extends Scheepers' Theorem \ref{schthm}
in the case that $X,Y\sbst\cC$.
\begin{thm}\label{MaddSMZ}
Assume that $X\in\Madd\cap\Pgal(\cC)$ and $Y\in\SMZ\cap\Pgal(\cC)$.
Then $X\x Y\in\SMZ$.
\end{thm}
\begin{proof}
By Theorem \ref{M&Nstarprod}, $\Madd \cap\Pgal(\cC)$ is \productive{}.
Recall that $\SMZ\cap\Pgal(\cC)$ is \productive{} too. By the
Galvin-Mycielski-Solovay Theorem (Proposition \ref{GMSk} for $k=1$),
the conditions of Theorem \ref{sum2prod} hold, and its consequence tells
what we are looking for.
\end{proof}

To prove the dual result, we need the following lemma.
For a set $J$ denote $J_x = \{y : (x,y)\in J\}$ and
$J^y = \{x : (x,y)\in J\}$.
Say that a family $\cJ$ which does not contain any $\cC^k$ as element
is a \emph{Fubini family} if,
whenever $J\in\cJ\cap\cC^{k+l}$, we have that
\begin{equation}\label{fub}
\{x\in\cC^k : J_x\nin\cJ\}\in\cJ\mbox{, and }\{y\in\cC^l : J^y\nin\cJ\}\in\cJ.
\end{equation}
The most important examples for Fubini families are
$\cM$ (Kuratowski-Ulam Theorem) and $\cN$ (Fubini Theorem).
To understand what we really prove, we
will say that $\cJ$ is a \emph{weakly Fubini family} if
``$\in\cJ$'' and ``$\in\cJ$'' in \eqref{fub} are replaced by
``$\neq\cC^k$'' and ``$\neq\cC^l$'', respectively.
Clearly, each Fubini family is a weakly Fubini family.

A set $X\sbst\cC^k$ is \emph{not $\cJ$-covering} if
for each $J\in\cJ\cap P(\cC^k)$, $X\oplus J\neq\cC^k$.

\begin{lem}
Assume that $\cJ$ is a weakly Fubini family.
Then the family of not $\cJ$-covering sets is \zeroproductive.
\end{lem}
\begin{proof}
Assume that $X\sbst\cC^k$ is not $\cJ$-covering,
$0\in\cC^l$, and $J\in\cJ\cap P(\cC^{k+l})$.
As $\cJ$ is a weakly Fubini family, there exists $y\in\cC^l$
such that $J^y\in\cJ$.
Thus, $\((X\x\{0\})\oplus J\)^y = X\oplus J^y\neq\cC^k$,
therefore $(X\x\{0\})\oplus J\neq\cC^{k+l}$.
\end{proof}

A set $X\sbst\cC^k$ is \emph{strongly meager} if
it is not $\cN$-covering.
Using the same proof as in Theorem \ref{MaddSMZ}, we get the following.
\begin{thm}
The product of a null-additive set in $\cC^k$ and a strongly meager
set in $\cC^l$ is strongly meager.\hfill\qed
\end{thm}

It is folklore that a product of strong measure zero sets need not have strong measure zero
(e.g., \cite{MilSpec}), and that the product of strongly meager sets need not be strongly meager.
We could not find a reference in the literature for the latter fact.
To see that it holds, we say that a set $S\sbst\cC$ is \emph{$\kappa$-Sierpi\'nski}
if $|S|\geq\kappa$ but for each null set $N$, $|S\cap N|<\kappa$.
Observe that the diagonal is null in $\cC^2$.

\begin{thm}
Assume that $\cov(\Null)=\b=\c$.
Then there exists a strongly meager set of reals $S\sbst\cC$
such that $S^2\oplus\Delta=\cC^2$,
where $\Delta=\{(x,x) : x\in\cC\}$.
\end{thm}
\begin{proof}
Since $\cov(\Null)=\c$ we can construct,
as in \cite[Lemma 42]{CBC}, a $\cov(\Null)$-Sierpi\'nski set $S$ such that $S\oplus S=\cC$.
Then $S^2\oplus\Delta=\cC^2$:
Given $y,z\in \cC$, choose $s,t\in S$ such that $s\oplus t=y\oplus z$,
and take $x=s\oplus y$. Then $s\oplus x=y$, and
$$t\oplus x=t\oplus (s\oplus s)\oplus x=(s\oplus t)\oplus (s\oplus x)=(y\oplus z)\oplus y=z,$$
thus $(s,t)\oplus(x,x)=(y,z)$.
Since $\b=\cov(\Null)$, $S$ is a $\b$-Sierpi\'nski set and
by \cite[p.~376]{CBC}, every Borel image of $S$ is bounded. Moreover, for each null set
$N$, $|S\cap N|<\cov(\Null)$. By a result of Pawlikowski
(see \cite{jubar} -- Definition 8.5.7, the observation after it, and Theorem 8.5.12),
these two properties imply that $S$ is strongly meager.
\end{proof}

\section{Products in the Euclidean space}\label{euclidprods}

As the mapping from $\R^k\x\R^k$ to $\R^k$ defined by
$(x,y)\mapsto x+y$ is Lipschitz,
Theorem \ref{prod2sum} remains true in the \emph{Euclidean space} $\R^k$.
However, we are unable to prove Theorem \ref{sum2prod} (in its current form)
for the Euclidean space ($\<\R^k,+\>$ and $\<\R,+\>$ are not
homeomorphic: $\R^k$ remains connected after removing a point).
We can, though, obtain similar results.

\begin{abbr}
A collection $\cP\sbst\Pgal(\R)$
is \emph{bi-\zeroproductive{}} if
for each $k$, $l$, and $X\in\cP\cap\R^k$,
if $0$ is the zero element of $\R^l$, then $X\x\{0\},\{0\}\x X\in\cP$.
\end{abbr}

\begin{lem}\label{sum2prodR}
Assume that $\cI,\cJ,\cK\sbst\Pgal(\R)$ and $\cI,\cJ$ are bi-\zeroproductive{}.
If
\bq
for each $k$ and all $X\in\cI\cap P(\R^k)$ and $Y\in\cJ\cap P(\R^k)$, $X+Y\in\cK$,
\eq
then
\bq
for each $X\in\cI$ and $Y\in\cJ$, $X\x Y,Y\x X\in\cK$.
\eq
\end{lem}
\begin{proof}
Assume that $X\in\cI\cap P(\R^k)$ and $Y\in\cJ\cap P(\R^l)$.
Then
$$X\x Y = (X\x\{0\})+(\{0\}\x Y).$$
As $\cI$ and $\cJ$ are bi-\zeroproductive{},
$X\x\{0\}\in\cI\cap P(\R^{k+l})$ and $\{0\}\x Y\in\cJ\cap P(\R^{k+l})$,
therefore $X\x Y\in\cK$. Similarly, $Y\x X\in\cK$.
\end{proof}

\begin{thm}\label{prods*R}
Assume that $\cI,\cJ,\cK\sbst\Pgal(\R)$, and that
$(\cJ,\cK)^\star$ and $(\cI,\cJ)^\star$ are bi-\zeroproductive{}.
Then for each $X\in(\cJ,\cK)^\star$ and $Y\in(\cI,\cJ)^\star$, $X\x Y,Y\x X\in (\cI,\cK)^\star$.
In particular, if $\cI^\star$ is bi-\zeroproductive{}, then
$\cI^\star$ is closed under taking finite products.
\end{thm}
\begin{proof}
Assume that $X\in(\cJ,\cK)^\star\cap P(\R^k)$, $Y\in(\cI,\cJ)^\star\cap P(\R^k)$,
and $I\in\cI\cap P(\R^k)$.
Then $Y+I\in\cJ$ and therefore $X+(Y+I)\in\cK$.
Thus $X+Y\in(\cI,\cK)^\star$.
As $(\cJ,\cK)^\star$ and $(\cI,\cJ)^\star$ are bi-\zeroproductive{},
our theorem follows from Lemma \ref{sum2prodR}.
\end{proof}

It was noticed by Marcin Kysiak that one can use our arguments
to obtain Scheepers' Theorem \ref{schthm} in the Euclidean space.
To see this, one simply has to generalize the corresponding theorem on
sums \cite{NSW} from $k=1$ to arbitrary $k$. The generalization is straightforward.
In Section \ref{completekysiac} we show that in fact, this
generalization is not necessary.

\medskip

The analogue of Theorem \ref{M&Nstarprod} in the Euclidean space
does not seem to follow from the results in this paper.

\begin{prob}
Assume that $X\sbst\R$ is meager- (respectively, null-) additive.
Does it follow that $X\x\{0\}$ is meager- (respectively, null-) additive?
\end{prob}

\begin{rem}
We can prove that every meager-additive subset of the Cantor
space, when viewed as a subset of $\R$, is meager-additive (with
respect to the usual addition in $\R$); and similarly for
null-additive (in both cases, the other direction is still open).
Consequently, the classes of meager-additive and null-additive
each contains a nontrivial subclass which is preserved under
taking finite products. We plan to treat this result elsewhere.
\end{rem}

\section{The Euclidean space through the looking glass}\label{completekysiac}

The results in Section \ref{euclidprods} are not easy to use,
as one should verify first that the additive
results given in the literature for $\R$ actually hold
in $\R^k$ for all $k$.
We suggest here another approach, which covers some of the
cases of interest.

\begin{definition}\label{binary}
The function $\Binary:\Cantor\to [0,1]$ is defined by
$$x \mapsto \sum_{i\in\N}\frac{x(i)}{2^{i+1}} = 0.x(0)x(1)x(2)\dots,$$
where the last term is in base $2$.
Let $C$ denote the collection of all eventually constant
elements of $\Cantor$, and $\Q_2=\Binary[C]$ denote the $2$-adic
rational numbers in $[0,1]$.
\end{definition}

\begin{lem}[folklore]\label{basicPhi}
~\be
\i $\Binary$ is a uniformly continuous surjection.
\i $C$ is countable.
\i $\Binary:\Cantor\sm C\to [0,1]\sm\Q_2$ is a homeomorphism
which preserves measure in both directions.
\ee
\end{lem}
\begin{proof}
(1) $\Binary$ is continuous on its compact domain, and clearly it is
onto.

(2) is obvious, and the only nontrivial part of (3) is that
$\Binary\inv$ is continuous on $[0,1]\sm\Q_2$ (it is not
\emph{uniformly} continuous: Take $x_n=0.10^n\closure{01}$ and
$y_n=0.01^n\closure{01}$, then $x_n-y_n\to 0$, but $d(\Binary\inv(x_n),\Binary\inv(y_n))=1$ for all $n$.).
Write $\tilde y$ for $\Binary\inv(y)$.
If $x_n\to x$ are elements of $[0,1]\sm\Q_2$ then from some
$n$ onwards, the $\tilde x_n(0)=\tilde x(0)$:
Assume that this is not the case. Then by moving to a
subsequence we may assume that
for all $n$, $\tilde x(0)\neq \tilde x_n(0)$.
Assume that $\tilde x(0)=0$ and $\tilde x_n(0)=1$ for all $n$
(the other case is similar).
Let $k=\min\{m>0 : x(m)=0\}$ (recall that $\tilde x$ is not eventually constant).
Then $x=0.01^{k-1}0\dots$, thus
$$x_n-x \ge 0.1-0.01^{k-1}0\closure{1} = 0.1-0.01^{k}=0.0^{k}1=\frac{1}{2^{k+1}}$$
for all $n$, a contradiction.

An inductive
argument shows that for each $k$, $\tilde x_n(k)=\tilde x(k)$
for all large enough $n$.
\end{proof}

In principle, we can use the function $\Binary$ to translate questions
about products in $\R$ into questions about products in $\cC$,
apply the results of Section \ref{cantorprods}, and translate back
to $\R$. The problem is that in our test-case,
Scheepers' Theorem \ref{schthm} in $\R$, we must deal with strong measure zero sets.
By (1) of Lemma \ref{basicPhi}, if $Y$ has strong measure zero then so does
$\Binary[Y]$. The other direction does not follow from Lemma
\ref{basicPhi}, since a homeomorphic image of a strong measure zero
set need not have strong measure zero \cite{ROTH41} (see \cite{ict} for a simple proof).

\begin{prop}\label{smzPhi}
For each $k$, $\Binary^k:\cC^k\to [0,1]^k$ preserves strong measure zero sets in both directions.
\end{prop}
\begin{proof}
Since $\Binary$ is uniformly continuous,
$\Binary^k$ is uniformly continuous. Thus, if
$X\sbst\cC^k$ has strong measure zero, then
so does $\Binary^k[X]$.

We now prove the other direction. The following is an
easy exercise.
\begin{lem}[folklore]\label{fewepsilons}
Assume that there exists $f:\N\to\N$ such that the metric space $(X,d)$
has the following property: For each sequence $\seq{\epsilon_n}$ of positive
reals, there exist a cover $\{I^n_m : n\in\N, m\le f(n)\}$ of $X$
satisfying $\diam(I^n_m)<\epsilon_n$ for each $n$ and $m$.
Then $(X,d)$ has strong measure zero.\hfill\qed
\end{lem}

To use Lemma \ref{fewepsilons}, we make the following observation.

\begin{lem}\label{fewepsilons2}
Assume that $I\sbst [0,1]$ and $\diam(I)<\epsilon$. Then
there exist $A_0,A_1\sbst\cC$, both with diameter $<\epsilon$,
such that $\Binary\inv[I]\sbst A_0\cup A_1$.
\end{lem}
\begin{proof}
Let $n$ be the maximal such that $2^{n-1}\epsilon\le 1$.
Then there exists $k<2^n-1$ such that
$$I\sbst \left [\frac{k}{2^n},\frac{k+2}{2^n}\right ]=
\left [\frac{k}{2^n},\frac{k+1}{2^n}\right ]\cup \left [\frac{k+1}{2^n},\frac{k+2}{2^n}\right ].$$
There exist sequences $s_0,s_1\in\two^n$ such that for $i=0,1$,
if $A_i=\{x\in\cC : s_i\sbst x\}$, then
$\Binary[A_i]=[(k+i)/2^n,(k+i+1)/2^n]$.
\end{proof}
Assume that $Y\sbst [0,1]^k$ has strong measure zero,
and $\seq{\epsilon_n}$ is a sequence of positive reals.
Let $I_n\sbst [0,1]^k$ be such that $\diam(I_n)<\epsilon_n$
and $Y\sbst\Union_nI_n$.
Fix $n$. Then $I_n$ is contained
in the product of $k$ intervals, $I^n_1,\dots I^n_k\sbst [0,1]$,
each with diameter $<\epsilon_n$.
For each $i=1,\dots,k$, use Lemma \ref{fewepsilons2}
to obtain sets $A^{n,i}_0,A^{n,i}_1$ as in the lemma.
Then
$$I_n\sbst\Union_{s\in\two^k}\Pi_{i=1}^k\Binary[A^{n,i}_{s(i)}]=
\Binary^k\left [\Union_{s\in\two^k}\Pi_{i=1}^kA^{n,i}_{s(i)}\right ],$$
so that $\Binary^{-k}[I_n]$ is covered by $2^k$ sets of diameter $<\epsilon_n$.
Since the sets $I_n$ cover $Y$, we have by Lemma
\ref{fewepsilons} that $\Binary^{-k}[Y]$ has strong measure zero.
\end{proof}

We now show how to prove Scheepers' Theorem in the Euclidean space.
\begin{lem}\label{avoidctbl}
Assume that $\cI$ is preserved under taking closed subsets, uniformly continuous images,
and countable unions, and that $\R\nin\cI$.
Then for each $X\in\cI\cap P(\R^k)$ and a countable set $Q\sbst\R$ there exists $x\in \R^k$
such that $(x+X)\cap \Union_{m<k}\R^m\x Q\x\R^{k-m-1}=\emptyset$.
\end{lem}
\begin{proof}
The assumptions imply that for each $i<k$, the projection $X_i = \pi_i[X]$ on the $i$th coordinate
is a member of $\cI$. As $Q$ is countable, $Q-X_i\neq\R$. Choose
$x_i\nin Q-X_i$. Then $(x_i+X_i)\cap Q=\emptyset$. Take
$x =(x_0,\dots,x_{k-1})$. Then $(x+X)\cap \Union_{m<k}\R^m\x Q\x\R^{k-m-1}=\emptyset$.
\end{proof}
Assume that $X\in\HURE\cap\SMZ\cap P(\R^k)$ and $Y\in\SMZ\cap P(\R^l)$.
It is well known (and easy to see)
that $\HURE$ and $\SMZ$ satisfy the assumptions of Lemma \ref{avoidctbl}.
Take
$$Q=\Union_{m,n\in\Z} (m\cdot \Q_2+n).$$
Then by Lemma \ref{avoidctbl}, we may assume that
$X$ is disjoint from $\Union_{m<k}\R^m\x Q\x\R^{k-m-1}$, and
$Y$ is disjoint from $\Union_{m<l}\R^m\x Q\x\R^{l-m-1}$.

For each $n$, set $X_n = X\cap [-n,n]^k$ and $Y_n = Y\cap [-n,n]^l$.
Then $X_n$ is a closed subset of $X$ and therefore has
the Hurewicz property. Moreover, $X_n$ and $Y_n$ have strong measure zero,
and $X\x Y=\Union_n X_n\x Y_n$.
Since $\SMZ$ is preserved under countable unions,
it is enough to show that $X_n\x Y_n$ has strong measure zero
for each $n$. Transforming (each coordinate of) $X_n,Y_n$ with the bi-Lipschitz homeomorphism
$x\mapsto (n+x)/2n$, we may assume (by our choice of $Q$!) that
$X_n, Y_n\sbst ([0,1]\sm\Q_2)^k$.

So assume that $X\in\HURE\cap\SMZ\cap P(([0,1]\sm\Q_2)^k)$, and $Y\in\SMZ\cap P(([0,1]\sm\Q_2)^l)$
Since $\Binary^k:\cC^k\to ([0,1]\sm\Q_2)^k$ is a homeomorphism,
$\Binary^{-k}[X]$ has the Hurewicz property.
By Proposition \ref{smzPhi}, $\Binary^{-k}[X],\Binary^{-l}[Y]\sbst\cC$ have strong measure zero.
By Scheepers' Theorem in $\cC$, $\Binary^{-k}[X]\x\Binary^{-l}[Y]$ has strong measure zero.
As $\Binary^k\x\Binary^l$ is uniformly continuous, $X\x Y$ has strong measure zero.

\section{Borel's Conjecture and conjunction of properties}\label{BC}

To put things in a wider context, we briefly describe the general framework.
Let $X$ be a topological space.
An open cover $\cU$ of $X$ is
an \emph{$\omega$-cover} of $X$ if $X$ is not in $\cU$ and for
       each finite subset $F$ of $X$, there is
       a set $U\in\cU$ such that $F\subseteq U$.
$\cU$ is a \emph{$\gamma$-cover} of $X$ if it is infinite and for each $x$ in
       $X$, $x\in U$ for all but finitely many $U\in\cU$.
Let $\O$, $\Omega$, and $\Gamma$ denote the collections of all countable open
covers, $\omega$-covers, and $\gamma$-covers of $X$, respectively.
Let $\scrA$ and $\scrB$ be collections of covers of a space $X$.
Following are selection hypotheses which $X$ might satisfy or not
satisfy.
\begin{itemize}
\item[$\sone(\scrA,\scrB)$\index{$\sone(\scrA,\scrB)$}:]
For each sequence $\seq{\cU_n}$ of members of $\scrA$,
there exist members $U_n\in\cU_n$, $n\in\N$, such that $\seq{U_n}\in\scrB$.
\item[$\sfin(\scrA,\scrB)$\index{$\sfin(\scrA,\scrB)$}:]
For each sequence $\seq{\cU_n}$
of members of $\scrA$, there exist finite (possibly empty)
subsets $\cF_n\sbst\cU_n$, $n\in\N$, such that $\Union_{n\in\N}\cF_n\in\scrB$.
\item[$\ufin(\scrA,\scrB)$\index{$\ufin(\scrA,\scrB)$}:]
For each sequence $\seq{\cU_n}$ of members of $\scrA$
\emph{which do not contain a finite subcover},
there exist finite (possibly empty) subsets $\cF_n\sbst\cU_n$, $n\in\N$,
such that $\seq{\cup\cF_n}\in\scrB$.
\end{itemize}

$\ufin(\O,\Gamma)$ is the Hurewicz property,
$\sfin(\O, \O)$ is the Menger property,
$\sone(\O, \O)$ is Rothberger's property $C''$,
and $\sone(\Omega,\Gamma)$ is the $\gamma$-property.

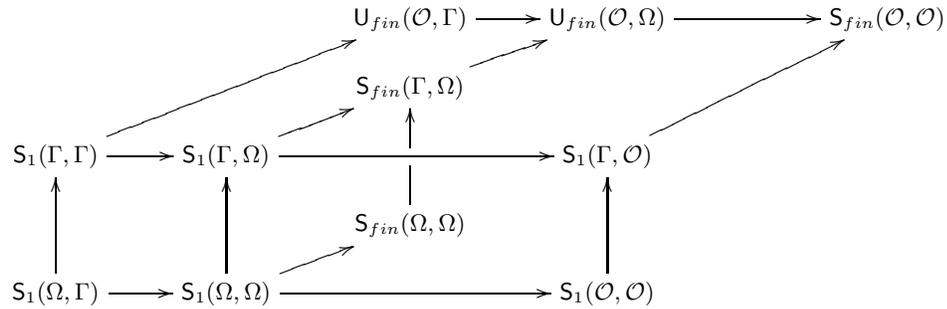
\begin{figure}[!ht]
{\scriptsize
\begin{center}
$\xymatrix@R=10pt{
&
&
& \ufin(\O,\Gamma)\ar[r]
& \ufin(\O,\Omega)\ar[rr]
& & \sfin(\O,\O)
\\
&
&
& \sfin(\Gamma,\Omega)\ar[ur]
\\
& \sone(\Gamma,\Gamma)\ar[r]\ar[uurr]
& \sone(\Gamma,\Omega)\ar[rr]\ar[ur]
& & \sone(\Gamma,\O)\ar[uurr]
\\
&
&
& \sfin(\Omega,\Omega)\ar'[u][uu]
\\
& \sone(\Omega,\Gamma)\ar[r]\ar[uu]
& \sone(\Omega,\Omega)\ar[uu]\ar[rr]\ar[ur]
& & \sone(\O,\O)\ar[uu]
}$
\caption{The Scheepers Diagram}\label{SchDiagram}
\end{center}
}
\end{figure}

Many equivalences hold among these properties, and the surviving ones
appear in Figure \ref{SchDiagram} (where an arrow denotes implication),
to which no arrow can be added except perhaps from
$\ufin(\O,\Gamma)$ or $\ufin(\O,\Omega)$ to $\sfin(\Gamma,\Omega)$
\cite{coc2}.

Let us write $BC(P)$ for the Borel Conjecture for \emph{metric
spaces} with property $P$, that is, the hypothesis that every
metric space with property $P$ is countable. Laver proved that
$BC(\SMZ)$ is consistent. Since $\sone(\O,\O)$ implies strong
measure zero, it follows that $BC(\sone(\Omega,\Gamma))$,
$BC(\sone(\Omega,\Omega))$, and $BC(\sone(\O,\O))$ are all
consistent. On the other hand, all other classes in the Scheepers
Diagram provably contain uncountable sets of reals \cite{coc2, ideals},
and therefore cannot satisfy $BC$.

In \cite{MillerBC} Miller proves that
$BC(\sone(\O,\O))$ implies
$BC(\SMZ)$, but $BC(\sone(\Omega,\Gamma))$ does not imply $BC(\SMZ)$.
We will extend this result.
With regards to the Scheepers Diagram \ref{SchDiagram}, the
best we can get is that
$BC(\sone(\Omega,\Omega))$ implies $BC(\SMZ)$
(and therefore $BC(\sone(\Omega,\Gamma))$ does not imply $BC(\sone(\Omega,\Omega))$\,).
We will prove a stronger result.

\begin{definition}
For a fixed topological space $X$,
$\Omega^\gpbl$ denotes the collection of open $\omega$-covers $\cU$ of $X$ such that:
There exists a partition $\cP$ of $\cU$ into finite sets such that
for each finite $F\sbst X$ and all but finitely many $\cF\in\cP$,
there exists $U\in\cF$ such that $F\sbst U$ \cite{coc7}.
\end{definition}

$\sone(\Omega,\Omega^\gpbl)$ is strictly stronger than $\sone(\Omega,\Omega)$ \cite{coc7}.

\begin{thm}\label{weakBC}
$BC(\sone(\Omega,\Omega^\gpbl))$ implies (and is therefore equivalent to) $BC(\SMZ)$.
\end{thm}
\begin{proof}
If $\aleph_1=\b$ then by \cite{ideals}
there exists an uncountable set of reals $X$ satisfying $\sone(\Omega,\Omega^\gpbl)$.

Assume that $\aleph_1<\b$, and $BC(\SMZ)$ fails.
Take any strong measure zero set $X$ with $|X|=\aleph_1$.
Then $|X|<\c$ and by a result of Carlson \cite[Lemma 8.1.9]{jubar},
we may assume that $X\sbst\R$.
As $|X|<\b$, all finite powers of $X$ have the Hurewicz property.
By Scheepers' Theorem \ref{schthm}, $X^2=X\x X$ has strong measure zero,
therefore $X^3=X\x X^2$ has strong measure zero, etc.

By \cite{coc7},
$\sone(\Omega,\Omega^\gpbl)$ is equivalent to all finite powers
having strong measure zero and satisfying the Hurewicz property.
\end{proof}

The arguments in the last proof actually establish the following.
\begin{thm}
For a set of reals $X$, the following are equivalent:
\be
\i\label{smz}
$X$ satisfies $\sfin(\Omega,\Omega^\gpbl)$ and has strong measure-zero,
\i\label{roth}
$X$ satisfies $\sfin(\Omega,\Omega^\gpbl)$ and $\sone(\O,\O)$,
\i\label{madd}
$X$ satisfies $\sfin(\Omega,\Omega^\gpbl)$ and is meager-additive,
\i\label{sakai}
$X$ satisfies $\sfin(\Omega,\Omega^\gpbl)$ and $\sone(\Omega,\Omega)$,
\i\label{gnpow}
$X$ satisfies $\sone(\Omega,\Omega^\gpbl)$.
\ee
\end{thm}
\begin{proof}
Clearly, $\ref{gnpow}\Impl\ref{sakai}\Impl\ref{roth}\Impl\ref{smz}$,
and $\ref{madd}\Impl\ref{smz}$.

($\ref{smz}\Impl\ref{gnpow}$)
Assume that (\ref{smz}) holds.
In \cite{coc7} it is proved that
$\sfin(\Omega,\Omega^\gpbl)$ is equivalent to satisfying the
Hurewicz property $\ufin(\O,\Gamma)$ in all finite powers.
By Scheepers' Theorem \ref{schthm}, all finite powers of $X$ satisfy
$\ufin(\O,\Gamma)$ and have strong measure zero.
By \cite{coc7}, $X$ satisfies $\sone(\Omega,\Omega^\gpbl)$

($\ref{smz}\Impl\ref{madd}$)
In \cite{NSW} it is proved that every strong measure zero set of reals with the
Hurewicz property is meager additive.
\end{proof}
The theorem also holds when $X\sbst\cC$.
In this case, the quoted assertion in the last proof can be proved directly
-- see Theorem \ref{(*)impliesM*} in the appendix.

\subsubsection*{Acknowledgements}
We thank Marcin Kysiak for his comment following Theorem
\ref{prods*R}, and
Tomek Bartoszy\'nski for the permission to include here
his proof of Theorem \ref{barthm}.

\appendix

\section{Direct proofs of quoted theorems}

Following is Bartoszy\'nski's (unpublished) combinatorial proof of
Scheepers' Theorem \ref{schthm} in $\cC$.
The proof uses the following characterization of strong measure zero.

\begin{lem}[{\cite[Lemma 8.1.13]{jubar}}]\label{smzchar}
For $X\sbst\cC$, the following are equivalent:
\be
\i $X$ has strong measure zero,
\i For each $f\in\NN$ there exists a function $g$
such that $g(n)\in {\two^{f(n)}}$ for all $n$, and
for each $x\in X$ there exist infinitely many $n$ such that
$x\restriction f(n)=g(n)$,
\i For each increasing sequence $\seq{m_n}$ there exists $z\in\cC$
such that for each $x\in X$ there exist infinitely many $n$ such that
$x\rest [m_n,m_{n+1}) = z\rest [m_n,m_{n+1})$.
\ee
\end{lem}

Let $\NNup$ denote the subspace of the Baire space $\NN$ consisting of
the increasing functions in $\NN$.
In \cite{hureslaloms} it is proved that if $X$ has the Hurewicz property
and $\Psi:X\to\NNup$ is continuous, then $\Psi[X]$ admits some slalom $h\in\NNup$,
that is, such that for each $x\in X$ and all but finitely many $n$, there exists $k$
such that $h(n)\le\Psi(x)(k)<h(n+1)$. This fact will be used in the
proof.

\begin{thm}\label{barthm}
Assume that $X\sbst\cC$ has the Hurewicz property and
strong measure zero, and $Y\sbst\cC$ has strong measure zero.
Then $X\x Y$ has strong measure zero.
\end{thm}
\begin{proof}[Proof (Bartoszy\'nski)]
Fix $f\in\NNup$ and let $g$ be as in Lemma \ref{smzchar} for $X$ and $f$.
Define a function $\Psi:X\to\NNup$ so that for each $x\in X$,
$\Psi(x)$ is the increasing enumeration of the set $\{n : x \restriction f(n)=g(n)\}$.
Then $\Psi$ is continuous, thus there exists $h\in\NNup$ such that
for each $x\in X$ and all but finitely many $n$, there exists $k$
such that $h(n)\le\Psi(x)(k)<h(n+1)$.

Consider a mapping $\Phi$ defined on $Y$ by
$$\Phi(y)(n) = \<(g(k), y\restriction f(k)): h(n)\le k<h(n+1)\>.$$
Then $\Phi$ is uniformly continuous.
Thus (essentially, by Lemma \ref{smzchar})
there exists a function $r$ such that
for all $y\in Y$ there exist infinitely many $n$ such that $\Phi(y)(n)=r(n)$.
From $r$ we decode a function $s$ such that $s(n)\in {\two^{f(n)}}\x {\two^{f(n)}}$
by $s(k) = r(n)(k)$ where $n$ is such that $h(n)\le k<h(n+1)$.

Then for all $x \in X$ and $y\in Y$ there exist infinitely many $n$ such that
$(x\restriction f(n), y \restriction f(n))=s(n)$, which shows that $X\x Y$ has strong measure zero.
\end{proof}

Using the bi-Lipschitz transformations $\Psi_k$ of Section \ref{cantorprods}, we obtain
Scheepers' Theorem in $\Pgal(\cC)$ from Theorem \ref{barthm}.

We can prove a result which is stronger (in light of the previous sections).
Following is a direct, combinatorial proof of one of the main theorems
in \cite{NSW} when restricted to the Cantor space.

\begin{thm}[\cite{NSW}]\label{(*)impliesM*}
Assume that $X\in\Pgal(\cC)$, $X$ has the Hurewicz property $\ufin(\O,\Gamma)$,
and strong measure zero.
Then $X$ is meager-additive.
\end{thm}
\begin{proof}
By Section \ref{cantorprods}, it suffices to prove the result for
$X\sbst\cC$.

Assume that $\seq{m_n}$ is an arbitrary increasing sequence.
By Lemma \ref{M*char}, it suffices to
find a sequence $\seq{l_n}$ and $z\in\cC$ such that for each
$x\in X$ and all but finitely many $n$,
$l_n\le m_k<m_{k+1}\le l_{n+1}$ and $x\|[m_k,m_{k+1})=z\|[m_k,m_{k+1})$
for some $k$.

By Lemma \ref{smzchar}, there exists $z\in\cC$
such that for each $x\in X$ there exist infinitely many $n$ such that
$x\rest [m_n,m_{n+1}) = z\rest [m_n,m_{n+1})$.
Again, for each $x\in X$ let $\Psi(x)$ be the increasing enumeration of these $n$s,
and use the fact that $X$ has the Hurewicz property to find a slalom
$h\in\NNup$ for $\Psi[X]$.

Take $l_n=m_{h(n)}$ for each $n$.
Fix $x\in X$. Since $h$ is a slalom for $\Psi[X]$,
for all but finitely many $n$ there exists
$k$ such that for $i=\Psi(x)(k)$, $h(n)\le i<h(n+1)$.
Then
$$l_n=m_{h(n)}\le m_i<m_{i+1}\le m_{h(n+1)}=l_{n+1},$$
and by the definition of $\Psi(x)$, $x\| [m_i,m_{i+1})=z\|[m_i,m_{i+1})$.
\end{proof}

\end{document}